\documentclass[10pt,twoside,dvips]{article}
\usepackage[a5paper,
layoutwidth=148.5mm,
layoutheight=7.77817in,
top=10mm,
bottom=10mm,
inner=37pt,
outer=27pt,
footskip=4mm]{geometry}

\usepackage[utf8]{inputenc}
\usepackage[T1]{fontenc}
\usepackage{lmodern}
\usepackage{microtype}
\usepackage{fix-cm}
\usepackage{textcomp}

\usepackage[vskip=3.5pt]{quoting}

\usepackage{natbib}
\usepackage[german,english,french]{babel} 
\usepackage[french]{varioref}

\newcommand{\bibAnnote}[2]{%
  \begin{quoting}\noindent\slshape#2\end{quoting}}
\providecommand{\bibAnnoteFile}[1]{}

\usepackage{graphicx}

\usepackage[autolanguage]{numprint}

\usepackage{framed}
\newenvironment{Cadre}{\framed}{\endframed}

\usepackage{amsmath}
\usepackage{amsfonts}

\usepackage{indentfirst}
\PassOptionsToPackage{urlcolor=magenta,citecolor=blue,anchorcolor=blue,linkcolor=blue,menucolor=cyan}{hyperref}%
\usepackage[pdftitle=,pdfauthor=,pdfdisplaydoctitle=true,pdflang=fr-FR,bookmarksopen=false,breaklinks=true,unicode=true,plainpages=false,hyperindex=true,pdfstartview=FitH,colorlinks=true,pdfpagelabels=true,bookmarksnumbered=true]{hyperref}
\ifx\texorpdfstring\undefined\newcommand\texorpdfstring[2]{#1}\fi
\ifx\pdfbookmark[3]{}\fi

\usepackage[anythingbreaks]{breakurl}%
\providecommand\url[1]{\texttt{\detokenize{#1}}}

\newcommand\segment[1]{\vrule width.6pt height4pt depth0pt\kern-.6pt\vbox{\hrule width#1pt height2.3pt depth-1.7pt}\kern-.6pt\vrule width.6pt height4pt depth0pt\hskip.05em}
\newcommand\segmenti[1]{\vrule width.6pt height4pt depth0pt\kern-.6pt\vbox{\hrule width#1pt height2.3pt depth-1.7pt}}

\newcommand\unite{\kern.5pt\vrule width1pt depth.4pt height8pt\kern.5pt}

\usepackage{newunicodechar}

\newunicodechar{·}{\mycdot}
\newcommand{\mycdot}{\kern-.2em\textperiodcentered\kern-.2em\relax}

\title{Un cours \emph{Littérature et mathématiques}}
\author{Spangle Durac \and Claude Merker \and Stefan Neuwirth}
\date{}
\begin{document}

\maketitle\thispagestyle{empty}

Nous rendons compte d'une expérience d'enseignement menée de 2008 à
2012 à l'université de Franche-Comté, en réponse à l'injonction
ministérielle de proposer des unités transversales dans les maquettes
de diplôme, tout en sachant bien qu'elle désignait ainsi des cours
d'anglais et de techniques d'expression.

Le but de cette expérience a été de développer à la fois un discours
sur les mathématiques et un discours sur la littérature, deux discours
indépendants, mais chacun rempli de l'écart entre l'un des domaines et
l'autre, et attentif aux éclairs de la pensée qui jaillissent de l'un
vers l'autre, témoignages d'unité dans la démarche intellectuelle.

En voici un exemple: le rôle de la conscience, qui à chaque instant
nous permet de nous saisir du processus de notre enquête, et de faire
de ce processus l'objet de l'enquête. Elle est une condition de la
possibilité des mathématiques en tant que telles. C'est par elle que le
calcul devient un objet de l'algèbre ou de la logique. C'est la capacité
qu'elle a de se viser elle-même qui engendre les vertiges de \emph{La
  disparition} de Georges Perec. L'absence tonitruante dont traite ce
roman renvoie aussi à l'inconscient. 

Ce compte rendu se veut un encouragement adressé aux amateur·e·s et un
éloge de l'amateurisme. Nous avons beaucoup réfléchi aux rapports
entre littérature et mathématiques, mais aucun·e de nous n'a le double
titre qui légitimerait ce cours selon le canon universitaire. Nous
revendiquons la gratuité de cette entreprise 
dans un contexte où la professionnalisation est le tombeau d'une
formation humaniste et holistique. Par le souci de ne pas
instrumentaliser les textes, littéraires et mathématiques, elle met
les étudiant·e·s dans une posture détendue, souple, attentive pour les
écouter.

Nous sommes convaincu·e·s de l'importance de faire un tel cours «sans
importance», pour l'un parce qu'il restaure l'unité de la culture,
pour l'autre parce qu'elle relève du jeu. Et nous nous sommes
pris·e·s au jeu au point d'organiser une séance de pâte à sel (voir
l'encadré sur le \emph{Traité de la roulette}).

\section{Légèreté, fragilité, fraicheur}
\label{sec:leger-frag-fraich}

Notre projet a été d'emblée d'introduire de la légèreté dans une
maquette de cours de mathématiques tous très sérieux et bien copieux,
une récréation plutôt qu'une charge de travail supplémentaire. La
légèreté est d'ailleurs une des valeurs retenues par Italo Calvino
pour le 21\ieme\ siècle dans ses \emph{Leçons américaines}! Nous avons
proposé aux étudiant·e·s d'engager une autre forme d'énergie que dans
les cours disciplinaires, faite de présence et d'attention pour
investir l'espace-temps du texte littéraire; nul travail n'a été
demandé en dehors des cours.

Nous sommes arrivé·e·s à la première séance très bien préparé·e·s par
de nombreuses réunions de travail dans lesquelles les
références littéraires et mathématiques ont fusé en réponse les unes aux autres, mais
notre parti a été d'emblée de tout oublier. Nous avons pris le risque de
séances fragiles, portées par des textes rendus à leur propre signification,
sans pointer vers un ailleurs. Nous avons parié qu'ils se
déploieraient mieux si nous les laissions s'exprimer librement.

Plutôt que de transmettre et de thésauriser des connaissances, notre
cours a pris la forme d'un happening qui se justifie par lui-même
et dans la mesure de sa justesse. Nos réunions préparatoires nous ont
permis de nous délester de tous les discours préalables, de tous les
prolégomènes, de toutes les conclusions; nous sommes retourné·e·s aux
fondements de ces discours, c'est-à-dire aux textes.

Les textes littéraires existent au moment de la lecture, et la
fragilité de notre dispositif
leur a évité le piédestal qui les aurait escamotés. Le
gain de cette démarche est la vie insufflée aux textes, par une
lecture neuve qui voit surgir des impressions fraiches.

\section{Le cours}
\label{sec:le-cours}

Les séances de cours ont été conçues en deux parties. Le premier tiers
est dédié à la mise en théâtre du \emph{Traité de la roulette} de
Blaise Pascal (voir l'encadré). Les deux autres tiers 
traitent d'un sujet,
d'une époque, d'un mouvement ou d'un·e auteur·e littéraire ou
mathématique.

Les deux premières années, nous avons évoqué la typographie (les
mathématiques du caractère et de la mise en page) et la métrique,
l'Oulipo (dont le jeu de la contrainte littéraire et le plaisir induit
de la combinatoire créent des conditions favorables pour l'apprenti
poète), \emph{Le cimetière marin} de Paul Valéry (muni·e·s du
commentaire de Gustave Cohen), les sonnets de Mallarmé (muni·e·s du
commentaire de Paul Bénichou, voir l'encadré), \emph{Sylvie et Bruno} de Lewis Carroll, \emph{Nous autres} de Ievgueni Zamiatine
et \emph{L'homme sans qualités} de Robert Musil.

Les protagonistes de ces deux derniers livres sont des
mathématiciens. Ils perçoivent le monde à travers leur formation d'ingénieur, et leurs mésaventures
questionnent l'adéquation du langage mathématique et de la rationalité
scientifique à la vie sociale et psychique.

Les deux années suivantes, nous avons travaillé notre matière à partir
de l'\emph{Intro\-duction à la littérature fantastique} de Tzvetan
Todorov. \emph{La Vénus d'Ille} de Prosper Mérimée, le \emph{Double
  assassinat dans la rue Morgue} d'Edgar Allan Poe et les
\emph{Chroniques martiennes} de Ray Bradbury nous ont permis de mener
une double enquête sur le mode d'existence des êtres de fiction
littéraires et mathématiques.

La littérature recourt aussi à la beauté froide des mathématiques pour décrire la place de l'individu face au monde: nommons à cet égard \emph{Le jeu des perles de verre} d'Hermann Hesse, \emph{Cristal de roche} d'Adalbert Stifter, \emph{Pierres} de Roger Caillois, \emph{La reine des neiges} de Hans Christian Andersen.

Par ailleurs, nous avons étudié des textes d'Euclide (voir l'encadré «Dire un» qui résulte d'une réflexion engagée après la fin de ce cours), d'Archimède, de Pierre de Fermat et de Johann Heinrich Lambert.

\section{Hamlet, Falstaff, Don Quichotte}

Évoquons un aspect qui rapproche littérature et mathématiques, et
constitue pour chacune un préalable à leur œuvre: c'est la capacité de
fiction. Cette notion est primordiale dans la théorie de la
littérature. Voici comment l'\emph{Introduction à la littérature fantastique} l'éclaire par une citation de Northrop Frye.
\begin{quoting}
  La littérature, comme les mathématiques, enfonce un coin
  dans l'antithèse de l'être et du non-être, qui est si importante
  pour la pensée discursive. [\dots] On ne peut dire de Hamlet et de
  Falstaff ni qu'ils existent ni qu'ils n'existent pas.
\end{quoting}
Complétons-la par la citation suivante de Roger Apéry.

\begin{quoting}
  Comme le platonicien et contrairement au formaliste, le
  mathématicien constructif reconnaît une certaine réalité aux objets
  mathématiques, mais les différencie essentiellement des objets
  matériels, en ne leur attribuant que les propriétés susceptibles de
  démonstration.  Une distinction analogue différencie les héros de
  roman des personnages historiques. Une question concernant
  Vercingétorix admet une réponse, même si elle échappe à nos moyens
  d’investigation; la même question concernant Don Quichotte n’a pas
  de réponse si celle-ci ne peut être déduite des affirmations du
  roman de Cervantes. En revanche, l’existence d’ensembles de réels
  plus nombreux que l’ensemble des entiers et moins nombreux que
  l’ensemble des réels n’a pas de réponse, car, comme Paul Cohen l’a
  démontré, ni cette existence ni sa négation ne peuvent être déduites
  des définitions usuelles des réels: l’ensemble des réels, comme Don
  Quichotte, est un être essentiellement incomplet.
\end{quoting}

L'enseignement des mathématiques arrive à éluder cette question en les
considérant comme une réalité ou comme une fiction dans leur ensemble,
dans laquelle on développe divers calculs dont on prétend à chaque
étape qu'ils ont lieu réellement au moment où ils sont effectués,
c'est-à-dire qu'ils correspondent à des faits. Mais il y a une étape à laquelle on paye le fait de ne pas avoir abordé cette question
de la fiction, c'est le raisonnement par récurrence. Dans l'encadré
«Dire un», nous présentons les nombres comme un jeu avec un seul signe, \unite, et deux
règles. Lorsqu'on veut savoir si quelque chose a lieu pour tous les
nombres, il faut «inverser» ces règles. Cela se fait ainsi: pour
qu'une propriété~$P$ vaille pour tous les nombres, il suffit de
prouver
\begin{itemize}
\item $P(\kern1pt\unite\kern1pt)$;
\item si $P(x)$, alors $P(x\kern1pt\unite\kern1pt)$.
\end{itemize}
Pour établir le deuxième point, il s'agit de faire une hypothèse,
$P(x)$, c'est-à-dire de changer consciemment l'état de son propre
entendement, faire \emph{comme si} $P(x)$, alors qu'on n'en sait rien,
et tenter d'aboutir dans cet état altéré à la
conclusion~$P(x\kern1pt\unite\kern1pt)$. Savoir faire cela, c'est
savoir feindre, \emph{i.\ e.}\ suspendre la réalité des faits,
\emph{i.\ e.}\ faire œuvre de fiction!
\[*\ *\ *\]

Ce compte rendu a d'abord voulu témoigner de notre enquête. Nous avons découvert que la littérature et les mathématiques ont le même rapport à la vérité, et que ce rapport est différent de celui des sciences de la nature et de la philosophe (voir l'encadré «Image et langage»).

L'expérience vécue a eu un impact notable sur le cours \emph{Musique
  et mathématiques}, qui a pris la relève depuis 2015. L'enjeu reste
de favoriser une perception sensible authentique et autonome par rapport à la
démarche intellectuelle. Mais nous avons appris entretemps qu'il faut
permettre aux étudiant·e·s de devenir activ·e·s en littérature, c'est-à-dire des
écrivain·e·s, pour que l'enquête soit profitable pour ce qui est de la création, littéraire et mathématique.

\emph{Nous remercions Alice Joly pour sa relecture.}

\renewcommand\bibpreamble{La bibliographie qui suit a pour but d'indiquer quelques-unes de nos trouvailles. Elles ne sont ni nécessaires ni suffisantes.}
\renewcommand\refname{Une sélection de références mathématico-littéraires}
\providecommand{\citenamefont}[1]{#1}
\providecommand{\urlprefix}{}
\providecommand{\selectlanguage}[1]{\relax}

\cite{*}

\clearpage

\begin{Cadre}
  \noindent\hfil\textbf{\large Image et langage}\bigskip

Lisons un texte littéraire tiré du début d'un
traité de Pascal Quignard sur Fronton, le maitre de Marc Aurèle, paru
dans le recueil \emph{Rhétorique spéculative}.

\begin{quoting}
  \noindent Fronton écrit à Marcus: «Il se trouve que le philosophe
  peut être imposteur et que l'amateur des lettres ne peut l'être. Le
  littéraire est chaque mot. D'autre part, son investigation propre
  est plus profonde à cause de l'image.» L'art des images --~que
  l'empereur Marc Aurèle nomme, en grec, icônes tandis que son maitre,
  Fronton, les nomme le plus souvent, en latin, images ou, à quelques
  reprises, en grec philosophique, métaphores~-- à la fois parvient à
  désassocier la convention dans chaque langue et permet de réassocier
  le langage au fond de la nature. [\dots] Fronton affirme que les
  arguments que peuvent articuler les philosophes ne sont que des
  claquements de langue parce qu'ils démontrent sans images: «S'il
  fait jour, alors il y a de la lumière.»
\end{quoting}

Ce texte parle de la littérature et de ce qui
l'oppose à la philosophie; les mathématicien·ne·s sont interpellé·e·s
et y reconnaissent leur discipline.

En effet, les mathématiques sont toujours prises dans l'œuvre de
«désassociation de la convention» et de «réassociation du langage au
fond de la nature», et sont guidées par les images. Le langage
mathématique est tout sauf stable, et il doit évoluer pour garder son
expressivité. Cela s'accompagne d'un réapprentissage permanent de la
signification des concepts, des formules, des diagrammes. On en
retient souvent un effet collatéral: les mathématiques forment une
tour de Babel dans laquelle les un·e·s ne comprennent les autres que
dans la mesure où ce travail de désassociation\slash réassociation est
répété, ou alors grâce aux images, aux métaphores.
\end{Cadre}

\clearpage

\begin{Cadre}
\noindent\hfil\textbf{\large Dire un}\bigskip
  
Lisons un tout petit texte mathématique: les deux premières
définitions du septième livre des \emph{Éléments} d'Euclide
d'Alexandrie.

\begin{quoting}
  1. Est \emph{unité} ce selon quoi chacune des choses existantes est
  dite une.

  2. Et un \emph{nombre} est la multitude composée d'unités.
\end{quoting}

On peut facilement rester insensible à ces définitions, tant la notion
de nombre semble claire et ces deux phrases obscures. Est-on alors
dans le cas d'un concept «premier», pour lequel toute tentative de
définition a pour effet de l'obscurcir? Nous sommes confronté·e·s au
premier mystère de la littérature: pourquoi écrire plutôt que se
taire et agir? L'usage du nombre est simple dans la mesure où on
l'emploie pour dénombrer et calculer, pourquoi alors vouloir faire du nombre
l'objet d'un discours?

En fait, il y a plusieurs raisons, mais retenons-en une seule: on a
besoin de narrer, et quand on narre une histoire, le moment critique
est le début, la première phrase qui instaure le rapport entre
auteur·e et auditeur·e.  La première définition ci-dessus a cet effet
fascinant de désigner, à travers la formulation «est dite», la personne
face aux choses comme étant celle qui effectue la définition de chaque chose
comme une unité. L'unité est le «dire un», le fait de reconnaitre une
chose existante comme étant une, de ne retenir de la
chose existante que l'unité qu'elle forme: c'est un acte de volonté,
et la définition de l'unité est réduite à cet acte. Il devient visible
lorsqu'on parle d'\emph{un} troupeau, d'\emph{une} myriade.

La deuxième définition est d'une autre nature: elle ne crée pas le
nombre, elle fournit seulement un mot pour la réitération du «dire
un». Mais la formulation entérine tacitement un constat: la multitude est une
multitude donnée et peu importe la manière dont on la compose ou
décompose. Il résulte d'une image, et on peut trouver la
trace de cette image dans les \emph{Éléments}: les nombres y sont
figurés par des segments, de sorte que trois nombres~\segment{20},
\segment{9}\ et~\segment{17}\ se juxtaposent en un seul
nombre~\segmenti{20}\segmenti{9}\segment{17}. Dans l'histoire qu'on
voudra narrer du nombre, on pourra employer le terme de tacite
\emph{associativité}, mais ce mot n'est attesté que depuis 1888
, suite à un questionnement sur
la notion même d'opération.

Multiplier un nombre, c'est le juxtaposer le nombre de fois voulu. Cette définition provient d'une image
de multitude: si un nombre est figuré par \segment9, alors son double
par \segmenti9\segment9\nolinebreak, son triple par \segmenti9\segmenti9\segment9,
son quadruple par \segmenti9\segmenti9\segmenti9\segment9, etc. La multiplication des
rapports ou des longueurs entre elles est un défi théorique
d'une autre nature.

Les définitions «Est \emph{unité} ce selon quoi chacune des choses
existantes est dite une», «Et un \emph{nombre} est la multitude
composée d'unités» nous font réfléchir au rapport complexe des
mathématiques avec la langue naturelle, française en
l'occurrence. Le concept de nombre peut être conçu comme un
jeu avec un seul signe, \unite, et deux règles:
\begin{itemize}
\item \unite;
\item si $x$, alors $x$\unite.
\end{itemize}
On construit \unite, et donc \unite\unite, et donc
\unite\unite\unite, et donc
\unite\unite\unite\unite, etc. Comment
capter la signification de ce jeu? On peut l'enraciner dans le langage
sous la forme d'un discours que l'on peut adresser à quelqu'un·e. Ce
discours épuise-t-il le jeu? Ou l'expérience de ce jeu est-elle
irréductible?

Notez bien que ces commentaires ne prétendent pas restituer la
signification que ces définitions ont eue dans l'antiquité
grecque. C'est-à-dire qu'en mathématiques comme en littérature, le
texte est vivant et mène sa vie propre au gré des lectures. Un des
charmes du retour régulier aux textes historiques est d'apercevoir
comment leur signification change du fait que nous changeons.

Ajoutons une dernière petite remarque sur les images. Les différents
discours sur le nombre s'accompagnent d'images différentes. Ces images
donnent une force et une prégnance au concept, et elles entretiennent
un rapport profond avec le texte: si celui-ci répond au défi de tenir
debout sans images, il rend néanmoins compte de leur pouvoir et
travaille à les susciter.
\end{Cadre}

\clearpage

\begin{Cadre}
\noindent\hfil\textbf{\large Le \emph{Traité de la roulette} de Blaise Pascal}\bigskip

Nous avons aussi lu, tout au long du quadrimestre, un texte
mathématique: le \emph{Traité de la roulette} de Blaise Pascal. Le
défi de cette lecture est de recréer la temporalité du texte, qui est
plutôt celle d'une pièce de théâtre. Dans ce texte, chaque mot
cherche sa place, non seulement dans le flot nécessaire des rapports
de cause à conséquence, mais aussi et surtout dans le cadre plus
général de la dramaturgie propre à l'œuvre: notre public doit
élaborer ce réseau de significations dans le temps réel de la
lecture.

Pour faire exister le \emph{Traité de la roulette}, pour le laisser opérer, il faut que le texte se présente de lui-même, il faut que ce soit Pascal qui parle, sans notre médiation, ou alors seulement par notre médiation matérielle. Pour rendre apparente la dimension théâtrale du texte, nous avons mis en place un dispositif:
\begin{itemize}
\item le texte a été lu et enregistré au préalable et diffusé;
\item le texte est aussi projeté sur un écran;
\item les figures géométriques et les formules sont tracées au tableau, au fur et à mesure qu'elles sont construites par le texte, dans un symbolisme cohérent avec la démarche de Pascal;
\item le tout avec le matériel le plus archaïque possible: magnétophone à bandes, rétroprojecteur, tableau noir et craie\dots
\end{itemize}
Nous avons néanmoins dû constater que la pédagogie, ainsi chassée,
revenait au galop, mue par notre pulsion professionnelle d'expliquer et de faire comprendre.

Nous avons particulièrement pris soin de maintenir le dialogue avec les
magnifiques figures géométriques de Pascal, conçues comme des machines que
l'imagination met en mouvement --~entre autres par l'usage de la pâte à sel!
Nous étions arrivé·e·s au passage un peu ardu des onglets et avions eu
besoin, pour le comprendre, de modeler les découpages que Pascal fait
dans son solide de révolution. Cet exercice pratique promettait d'être
profitable aux étudiant·e·s: nous leur avons donc apporté le matériau
adéquat. La pâte à sel a été utile à nos explications, mais son
introduction à l'université avait l'effronterie du gag, bienvenue à un
moment où nos fronts étaient devenus un peu trop studieux.

Voici comment Pascal introduit l'onglet dans un des avertissements de la \emph{Lettre de M. Dettonville à M. de Carcavi}:\nopagebreak

 \begin{quoting}
   On suppose ici toujours que le triligne est une figure plane, et
   que la courbe est de telle sorte que tant les sinus que les
   ordonnées ne la rencontrent qu'en un point. Et les portions de
   l'axe, de la base et de la courbe sont toutes égales tant entre
   elles que les unes aux autres.

   \noindent\hfil\includegraphics[width=50mm]{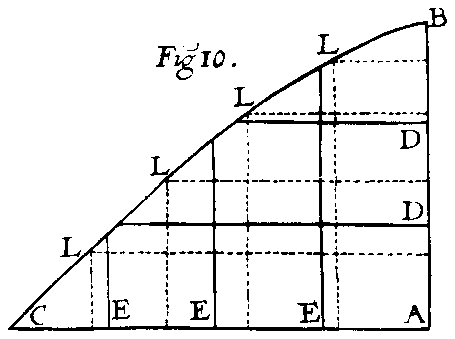}

   Il faut aussi remarquer que les sinus diffèrent des ordonnées, en
   ce que les sinus naissent des divisions égales de la courbe, et les
   ordonnées des divisions égales de l'axe ou de la base.

   Soient maintenant entendues des perpendiculaires élevées sur le
   plan de tous les points du triligne, qui forment un solide
   prismatique infini, qui aura le triligne pour base, lequel soit
   coupé par un plan incliné passant par l'axe ou par la base du
   triligne. La portion de ce solide retranchée par le plan
   s'appellera \emph{onglet}.
\end{quoting}

\end{Cadre}

\clearpage

\clearpage
\begin{Cadre}
  \begin{center}
    \large\bfseries Mallarmé
  \end{center}
  \label{sec:mallarme}

Quels points communs entre mathématiques et poésie? Leurs textes sont
repérables par leur forme:
\begin{itemize}
\item pieds, rimes, enjambements, inversions et syntaxe étrange, etc.\ pour la poésie.
\item formules, signes hors langue naturelle pour les mathématiques.
\end{itemize}

Poésie comme mathématiques ne sont pas tenues à la représentation (alors que la littérature représente, voir Todorov sur la littérature fantastique).

Un poème de Mallarmé, le sonnet en -yx, s'organise autour d'un mot, «ptyx», choisi \emph{a priori} pour une raison de rime, et non pour ce qu'il signifie, puisque Mallarmé a demandé à ses amis: «[\dots] je n'ai que trois rimes en ix, concertez-vous pour m'envoyer le sens réel du mot ptyx: on m'assure qu'il n'existe dans aucune langue, ce que je préfèrerais de beaucoup à fin de me donner le charme de le créer par la magie de la rime». On est loin du précepte de Boileau: «La rime est une esclave, et ne doit qu'obéir».

On s'attend à ce qu'un tel point de vue ait des répercussions sur le sens général du poème. Mallarmé, toujours parlant de son même sonnet: «il est inverse, je veux dire que le sens, s'il en a un (mais je me consolerais du contraire grâce à la dose de poésie qu'il renferme, ce me semble) est évoqué par un mirage interne des mots mêmes» et, plus loin, «j'ai pris ce sujet d'un sonnet nul et se réfléchissant de toutes les façons».

Deux résonances intérieures au poème sont bien visibles:
\begin{itemize}
\item Le mot «ptyx», directement importé du grec, est défectif, donc adapté à parler du néant, meilleur sens admissible pour ce sonnet «nul», \emph{i.\ e.}\ sans sens\dots
\item Le sonnet se termine par la réflexion du septuor (la Grande Ourse) dans un miroir, et justement les rimes des deux sortes, -or et -ix, sont chacune au nombre de sept.
\end{itemize}

Bernard Marchal, dans l'édition de la Pléiade, écrit: «Le sonnet en -yx est l'illustration parfaite d'une poésie désormais consciente d'elle-même qui consacre l'immanence du sens».

Le courant mathématique moderne du bourbakisme a de même mis l'accent sur le développement interne. Immanence plutôt que transcendance là aussi, avec tous les inconvénients afférents puisque les mathématiques, elles, ne peuvent raisonnablement se couper de la physique dont le rapport au réel est complètement nécessaire.
\end{Cadre}


\end{document}